\documentclass[11pt]{article}
\usepackage{bbm}
\usepackage{mathrsfs}
\usepackage{amscd}
\usepackage{amsmath,amsfonts,amssymb,amscd}
\usepackage{indentfirst,graphics,epsfig,psfrag}
\input{epsf}
\usepackage{ifpdf}
\usepackage{enumerate}
\usepackage{appendix}
\usepackage{enumerate}
\usepackage{color}
\usepackage{lineno}

\makeatletter \@addtoreset{figure}{section} \makeatother
\makeatletter
\long\def\@makecaption#1#2{%
   \vskip 10\p@
   \setbox\@tempboxa\hbox{{#1}\ \ #2}%
   \ifdim \wd\@tempboxa >\hsize

       {#1}\ \ #2\par
   \else
       \hbox to\hsize{\hfil\box\@tempboxa\hfil}%
   \fi}
\makeatother

\newtheorem{thm}{Theorem}
\newtheorem{cor}{Corollary}
\newtheorem{lem}{Lemma}

\newtheorem{obs}{Observation}
\newtheorem{pro}{Proposition}

\newcommand{\qed}{{\hfill\rule{3pt}{7pt}}}

\def\qed{\hfill \rule{4pt}{7pt}}

\begin{document}
\title{\textbf{The Steiner $4$-diameter of a graph} \footnote{Supported by the National Science Foundation of China
(Nos. 11601254, 11551001, 11371205, 11161037, 11101232, 11461054) and the Science Found of Qinghai Province
(Nos. 2016-ZJ-948Q, and 2014-ZJ-907).}}
\author{
\small Zhao Wang$^{1,2}$, \ \ Yaping Mao$^{2,4}$\footnote{Corresponding author}, \ \ Hengzhe Li$^{3}$, \ \ Chengfu Ye$^{2,4}$\\[0.2cm]
\small $^1$School of Mathematical Sciences, Beijing Normal\\
\small University, Beijing 100875, China\\[0.2cm]
\small $^2$Department of Mathematics, Qinghai Normal\\[-1mm]
\small University, Xining, Qinghai 810008, China\\[0.2cm]
\small $^3$School of Mathematical Sciences, Henan Normal\\
\small University, Xinxiang 453007, China\\[0.2cm]
\small $^4$Center for Mathematics and Interdisciplinary Sciences\\
\small  of Qinghai Province, Xining, Qinghai 810008, China\\[0.2cm]
\small E-mails: wangzhao@mail.bnu.edu.cn; maoyaping@ymail.com; \\
\small  hengzhe\_li@126.com; yechf@qhnu.edu.cn}
\date{}
\maketitle
\begin{abstract}
The Steiner distance of a graph, introduced by Chartrand,
Oellermann, Tian and Zou in 1989, is a natural generalization of the
concept of classical graph distance. For a connected graph $G$ of
order at least $2$ and $S\subseteq V(G)$, the \emph{Steiner
distance} $d_G(S)$ among the vertices of $S$ is the minimum size among
all connected subgraphs whose vertex sets contain $S$. Let $n,k$ be
two integers with $2\leq k\leq n$. Then the \emph{Steiner
$k$-eccentricity $e_k(v)$} of a vertex $v$ of $G$ is defined by
$e_k(v)=\max \{d(S)\,|\,S\subseteq V(G), \ |S|=k, \ and \ v\in S
\}$. Furthermore, the \emph{Steiner $k$-diameter} of $G$ is
$sdiam_k(G)=\max \{e_k(v)\,|\, v\in V(G)\}$. In 2011, Chartrand,
Okamoto and Zhang showed that $k-1\leq sdiam_k(G)\leq n-1$. In this
paper, graphs with $sdiam_4(G)=3,4,n-1$ are characterized,
respectively.\\[2mm]
\end{abstract}

\section{Introduction}

All graphs in this paper are undirected, finite and simple. We refer
to \cite{Bondy} for graph theoretical notation and terminology not
described here. For a graph $G$, let $V(G)$, $E(G)$,
$e(G)$, $\delta(G)$, and $\overline{G}$ denote the set of vertices, the set
of edges, the size, minimum degree, and the complement of $G$, respectively.
In this paper, we let $K_{n}$, $P_n$, $K_{1,n-1}$ and $C_n$ be the
complete graph of order $n$, the path of order $n$, the star of order $n$, and the
cycle of order $n$, respectively. For any subset $X$ of $V(G)$, let $G[X]$ denote the subgraph induced
by $X$; similarly, for any subset
$F$ of $E(G)$, let $G[F]$ denote the subgraph induced by $F$. We use
$G\setminus X$ to denote the subgraph of $G$ obtained by removing all the
vertices of $X$ together with the edges incident with them from $G$;
similarly, we use $G\setminus F$ to denote the subgraph of $G$
obtained by removing all the edges of $F$ from $G$. If $X=\{v\}$ and
$F=\{e\}$, we simply write $G-v$ and $G\setminus e$ for $G-\{v\}$
and $G\setminus \{e\}$, respectively. For two subsets $X$ and $Y$ of
$V(G)$ we denote by $E_G[X,Y]$ the set of edges of $G$ with one end
in $X$ and the other end in $Y$. If $X=\{x\}$, we simply write
$E_G[x,Y]$ for $E_G[\{x\},Y]$. We
divide our introduction into the following four subsections to state
the motivations and our results of this paper.

\subsection{Distance and its generalizations}

Distance is one of the most basic concepts of
graph-theoretic subjects. For a graph $G$, let $V(G)$, $E(G)$,
and $e(G)$ denote the set of vertices, the set
of edges, and the size of $G$, respectively. If $G$ is a connected graph and $u,v\in
V(G)$, then the \emph{distance} $d_G(u,v)$ between $u$ and $v$ is the
length of a shortest path connecting $u$ and $v$. If $v$ is a vertex
of a connected graph $G$, then the \emph{eccentricity} $e(v)$ of $v$
is defined by $e(v)=\max\{d_G(u,v)\,|\,u\in V(G)\}$. Furthermore, the
\emph{radius} $rad(G)$ and \emph{diameter} $diam(G)$ of $G$ are
defined by $rad(G)=\min\{e(v)\,|\,v\in V(G)\}$ and $diam(G)=\max
\{e(v)\,|\,v\in V(G)\}$. These last two concepts are related by the
inequalities $rad(G)\leq diam(G) \leq 2 rad(G)$. The \emph{center}
$C(G)$ of a connected graph $G$ is the subgraph induced by the
vertices $u$ of $G$ with $e(u)=rad(G)$. Recently, Goddard and
Oellermann gave a survey paper on this subject, see \cite{Goddard}.

The distance between two vertices $u$ and $v$ in a connected graph
$G$ also equals the minimum size of a connected subgraph of $G$
containing both $u$ and $v$. This observation suggests a
generalization of distance. The Steiner distance of a graph,
introduced by Chartrand, Oellermann, Tian and Zou in 1989, is a
natural and nice generalization of the concept of classical graph
distance. For a graph $G(V,E)$ and a set $S\subseteq V(G)$ of at
least two vertices, \emph{an $S$-Steiner tree} or \emph{a Steiner
tree connecting $S$} (or simply, \emph{an $S$-tree}) is a
subgraph $T(V',E')$ of $G$ that is a tree with $S\subseteq V'$. Let
$G$ be a connected graph of order at least $2$ and let $S$ be a
nonempty set of vertices of $G$. Then the \emph{Steiner distance}
$d_G(S)$ among the vertices of $S$ (or simply the distance of $S$)
is the minimum size among all connected subgraphs whose vertex sets
contain $S$. Note that if $H$ is a connected subgraph of $G$ such
that $S\subseteq V(H)$ and $|E(H)|=d_G(S)$, then $H$ is a tree.
Observe that $d_G(S)=\min\{e(T)\,|\,S\subseteq V(T)\}$, where $T$ is
subtree of $G$. Furthermore, if $S=\{u,v\}$, then $d_G(S)=d(u,v)$ is
the classical distance between $u$ and $v$. Set
$d_G(S)=\infty$ when there is no $S$-Steiner tree in $G$.

Let $n$ and $k$ be two integers with $2\leq k\leq n$. The
\emph{Steiner $k$-eccentricity $e_k(v)$} of a vertex $v$ of $G$ is
defined by $e_k(v)=\max \{d(S)\,|\,S\subseteq V(G), |S|=k,~and~v\in
S \}$. The \emph{Steiner $k$-radius} of $G$ is $srad_k(G)=\min \{
e_k(v)\,|\,v\in V(G)\}$, while the \emph{Steiner $k$-diameter} of
$G$ is $sdiam_k(G)=\max \{e_k(v)\,|\,v\in V(G)\}$. Note for every
connected graph $G$ that $e_2(v)=e(v)$ for all vertices $v$ of $G$
and that $srad_2(G)=rad(G)$ and $sdiam_2(G)=diam(G)$. Each vertex of
the graph $G$ of Figure 1 $(c)$ is labeled with its Steiner
$3$-eccentricity, so that $srad_3(G)=4$ and $sdiam_3(G)=6$.

\begin{obs}\label{obs2}
Let $k,n$ be two integers with $2\leq k\leq n$.

$(1)$ If $H$ is a spanning subgraph of $G$, then $sdiam_k(G)\leq
sdiam_k(H)$.

$(2)$ For a connected graph $G$, $sdiam_k(G)\leq sdiam_{k+1}(G)$.
\end{obs}

In \cite{ChartrandOZ}, Chartrand, Okamoto, Zhang obtained the
following result.
\begin{thm}{\upshape\cite{ChartrandOZ}}\label{th1}
Let $k,n$ be two integers with $2\leq k\leq n$, and let $G$ be a
connected graph of order $n$. Then $k-1\leq sdiam_k(G)\leq n-1$.
Moreover, the upper and lower bounds are sharp.
\end{thm}

In \cite{DankelmannSO2}, Dankelmann, Swart and Oellermann obtained a
bound on $sdiam_k(G)$ for a graph $G$ in terms of the order of $G$
and the minimum degree $\delta$ of $G$, that is, $sdiam_k(G)\leq
\frac{3n}{\delta+1}+3k$. Later, Ali, Dankelmann, Mukwembi
\cite{AliDM} improved the bound of $sdiam_k(G)$ and showed that
$sdiam_k(G)\leq \frac{3n}{\delta+1}+2k-5$ for all connected graphs
$G$. Moreover, they constructed graphs to show that the bounds are
asymptotically best possible.

As a generalization of the center of a graph, the \emph{Steiner
$k$-center} $C_k(G)\ (k\geq 2)$ of a connected graph $G$ is the
subgraph induced by the vertices $v$ of $G$ with $e_k(v)=srad_k(G)$.
Oellermann and Tian \cite{OellermannT} showed that every graph is
the $k$-center of some graph. In particular, they showed that the
$k$-center of a tree is a tree and those trees that are $k$-centers
of trees are characterized. The \emph{Steiner $k$-median} of $G$ is
the subgraph of $G$ induced by the vertices of $G$ of minimum
Steiner $k$-distance. For Steiner centers and Steiner medians, we
refer to \cite{Oellermann, Oellermann2, OellermannT}.

The \emph{average Steiner distance} $\mu_k(G)$ of a graph $G$,
introduced by Dankelmann, Oellermann and Swart in
\cite{DankelmannOS}, is defined as the average of the Steiner
distances of all $k$-subsets of $V(G)$, i.e.
$$
\mu_k(G)={n\choose k}^{-1}\sum_{S\subseteq V(G),|S|=k}d_G(S).
$$
For more details on average Steiner distance, we refer to
\cite{DankelmannOS, DankelmannSO}.

Let $G$ be a $k$-connected graph and $u$, $v$ be any pair of
vertices of $G$. Let $P_k(u,v)$ be a family of $k$ inner vertex-disjoint
paths between $u$ and $v$, i.e., $P_k(u,v)=\{P_1,P_2,\cdots,P_k\}$,
where $p_1\leq p_2\leq \cdots \leq p_k$ and $p_i$ denotes the number
of edges of path $P_i$. The \emph{$k$-distance} $d_k(u,v)$ between
vertices $u$ and $v$ is the minimum $p_k$ among all $P_k(u,v)$ and
the \emph{$k$-diameter} $d_k(G)$ of $G$ is defined as the maximum
$k$-distance $d_k(u,v)$ over all pairs $u,v$ of vertices of $G$. The
concept of $k$-diameter emerges rather naturally when one looks at
the performance of routing algorithms. Its applications to network
routing in distributed and parallel processing are studied and
discussed by various authors including Chung \cite{Chung}, Du, Lyuu
and Hsu \cite{Du}, Hsu \cite{Hsu, Hsu2}, Meyer and Pradhan
\cite{Meyer}.

\subsection{Application background of Steiner distance}

Let $G$ be a $k$-connected graph and $u$, $v$ be any pair of
vertices of $G$. Let $P_k(u,v)$ be a family of $k$ internally vertex-disjoint
paths between $u$ and $v$, i.e. $P_k(u,v)=\{P_{p_1},P_{p_2},\cdots,P_{p_k}\}$,
where $p_1\leq p_2\leq \cdots \leq p_k$ and $p_i$ denotes the number
of edges of path $P_{p_i}$. The \emph{$k$-distance} $d_k(u,v)$ between
vertices $u$ and $v$ is the minimum $|p_k|$ among all $P_k(u,v)$ and
the \emph{$k$-diameter} $d_k(G)$ of $G$ is defined as the maximum
$k$-distance $d_k(u,v)$ over all pairs $u,v$ of vertices of $G$. The
concept of $k$-diameter emerges rather naturally when one looks at
the performance of routing algorithms. Its applications to network
routing in distributed and parallel processing are studied and
discussed by various authors including Chung \cite{Chung}, Du, Lyuu
and Hsu \cite{Du}, Hsu \cite{Hsu, Hsu2}, Meyer and Pradhan
\cite{Meyer}.

The Wiener index $W(G)$ of the graph $G$ is defined as
$W(G)=\sum_{\{u,v\} \subseteq V(G)} d_G(u,v)$.
Details on this oldest distance--based topological index can be
found in numerous surveys, e.g., in
\cite{Rouv1, Rouv2, Dobrynin, Xu}.
Li et al. \cite{LMG} put forward a Steiner--distance--based generalization of the
Wiener index concept. According to \cite{LMG}, the {\it $k$-center Steiner Wiener
index\/} $SW_k(G)$ of the graph $G$ is defined by
\begin{equation}                    \label{sw}
SW_k(G)=\sum_{\overset{S\subseteq V(G)}{|S|=k}} d(S)\,.
\end{equation}
For $k=2$, the above defined Steiner Wiener index coincides with
the ordinary Wiener index. It is usual to consider $SW_k$ for
$2 \leq k \leq n-1$, but the above definition would be applicable
also in the cases $k=1$ and $k=n$, implying $SW_1(G)=0$ and
$SW_n(G)=n-1$. A chemical application of $SW_k$ was recently reported in \cite{GFL}.
Gutman \cite{GutmanSDD} offered an analogous
generalization of the concept of degree distance. Later, Furtula, Gutman, and Katani\'{c} \cite{FurtulaGK}
introduced the concept
of Steiner Harary index and gave its chemical applications. For more details on Steiner distance
indices, we refer to \cite{FurtulaGK, GFL, GutmanSDD, LMG, LMG2, MWG, MWGK, MWGL}.

\subsection{Our results}

From Theorem \ref{th1}, we have $k-1\leq sdiam_k(G)\leq n-1$. In \cite{Mao}, Mao characterized the graphs with
$sdiam_3(G)=2,3,n-1$, respectively, and studied the
Nordhaus-Gaddum-type problem of the parameter $sdiam_k(G)$.

In this paper, graphs with $sdiam_4(G)=3,4,n-1$ are characterized, respectively.
\begin{thm}\label{th2}
Let $G$ be a connected graph of order $n \ (n\geq 4)$.

$(i)$ If $n=4$, then $sdiam_4(G)=3$;

$(ii)$ If $n\geq 5$, then $sdiam_4(G)=3$ if
and only if $n-3\leq\delta(G)\leq n-1$ and $C_4$ is not a subgraph of $\overline{G}$.
\end{thm}

\begin{figure}[!hbpt]
\begin{center}
\includegraphics[scale=0.750]{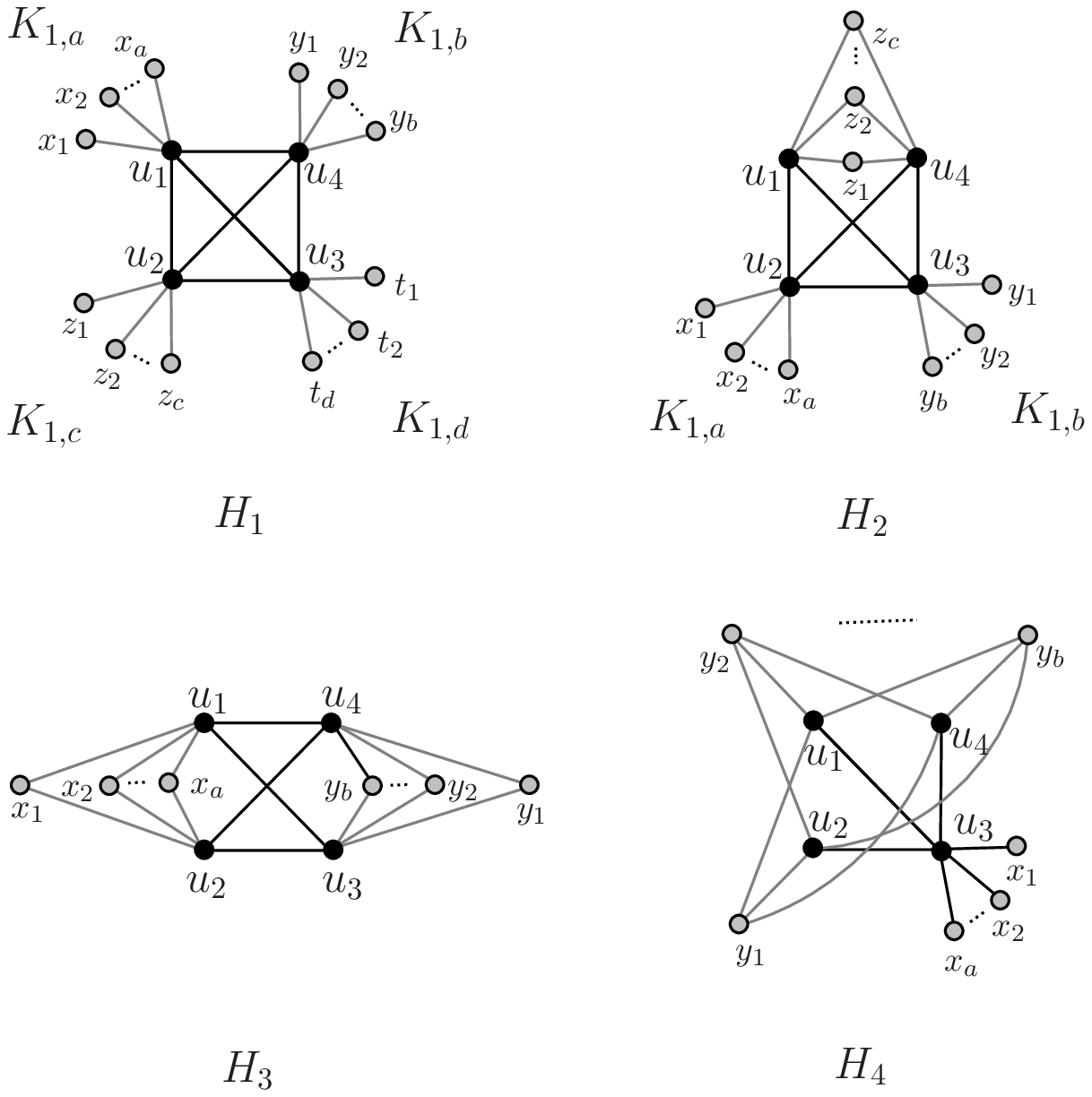}
\end{center}
\begin{center}
Figure 1: Graphs for Theorem \ref{th3}.\label{fig1}
\end{center}
\end{figure}

A graph $H_1$ is
defined as a connected graph of order $n \ (n\geq 5)$ obtained from a $K_4$ with vertex set $\{u_1,u_2,u_3,u_4\}$
and four stars $K_{1,a}, K_{1,b}, K_{1,c}, K_{1,d}$ by identifying the
center of one star and one vertex in $\{u_1,u_2,u_3,u_4\}$, where $0\leq a\leq
b\leq c\leq d$, $d\geq 1$, and $a+b+c+d=n-4$; see Figure \ref{fig1}.

A graph $H_2$ is
defined as a connected graph of order $n \ (n\geq 5)$ obtained from $K_4-e$ with vertex set
$\{u_1,u_2,u_3,u_4\}$, $e=u_1u_4$ and two stars $K_{1,a}$, $K_{1,b}$ by identifying the
center of a star and one vertex in $\{u_2,u_3\}$, and then adding the paths $u_1z_iu_4 \ (1\leq i\leq c)$, where $0\leq a\leq
b$, $b\geq 0$, $c\geq 0$ and $a+b+c=n-4$; see Figure \ref{fig2}.

A graph $H_3$ is
defined as a connected graph of order $n \ (n\geq 5)$ obtained from a cycle $C_4=u_1u_2u_3u_4u_1$ by adding the paths $u_1x_iu_2 \ (1\leq i\leq a)$ and the paths $u_3y_ju_4 \ (1\leq j\leq b)$, where $0\leq a\leq
b$, $b\geq 1$ and $a+b=n-4$; see Figure \ref{fig2}.

A graph $H_4$ is
defined as a connected graph of order $n \ (n\geq 5)$ obtained from a star $K_{1,3}$ with vertex set $\{u_1,u_2,u_3,u_4\}$
and a star $K_{1,a}$ by identifying $u_3$ and the center of $K_{1,a}$, where $u_3$ is the center of $K_{1,3}$, and then
adding the vertices $y_i$ and the edges $y_iu_j \ (1\leq i\leq b, \ j=1,2,4)$, where $0\leq a\leq
b$, $b\geq 1$ and $a+b=n-4$; see Figure \ref{fig2}.

\begin{thm}\label{th3}
Let $G$ be a connected graph of order $n \ (n\geq 5)$. Then $sdiam_4(G)=4$ if
and only if $G$ satisfies one of the following conditions.

$(i)$ $\delta(G)=n-3$ and $C_4$ is  a subgraph of $\overline{G}$;

$(ii)$ $\delta(G)\leq n-4$ and each $H_i\ (1\leq i\leq 4)$ is not a spanning subgraph of $\overline{G}$ (see Figure \ref{fig1}).
\end{thm}

We now define some graph classes.

\begin{itemize}
\item Let $T_{a,b,c,d} \ (0\leq a,b,c,d\leq n-1, a+b+c+d\leq n-1)$ be a tree of order $n  \ (n\geq 5)$ obtained from three paths
$P_1,P_2,P_3$ of length $n-b-c-1,b,c$ respectively by identifying
the $(a+1)$-th vertex of $P_1$ and one endvertex of $P_2$, and then
identifying the $(n-b-c-d)$-th vertex of $P_1$ and one endvertex of
$P_3$ (Note that $u$ and $v$ can be the same vertex);

\item Let $\triangle_{a,b,c,d} \ (0\leq a,b,c,d\leq n-2, a+b+c+d\leq n-2)$ be an unicyclic graph of order $n  \ (n\geq 5)$ obtained from three paths
$P_1,P_2,P_3$ of length $n-b-c-1,b+1,c$ respectively by identifying
the $(a+1)$-th vertex of $P_1$ and one endvertex of $P_2$, and then
identifying the $(n-b-c-d)$-th vertex of $P_1$ and one endvertex of
$P_3$, and then adding an edge $u_{b+1}v_{a+2}$  (Note that $v_{a+2}$ and $v$ can be the same vertex).

\item Let $\triangle_{a,b,c,d}' \ (0\leq a,b,c,d\leq n-3, \ a+b+c+d\leq n-3)$ be an bicyclic graph of order $n  \ (n\geq 5)$ obtained from three paths
$P_1,P_2,P_3$ of length $n-b-c-1,b+1,c+1$ respectively by identifying
the $(a+1)$-th vertex of $P_1$ and one endvertex of $P_2$, and then
identifying the $(n-b-c-d)$-th vertex of $P_1$ and one endvertex of
$P_3$, and then adding two edges $u_{b+1}v_{a+2}$ and $w_{c+1}x_{d+2}$ (Note that $v_{a+2}$ and $v$ can be the same vertex).

\item Let $G_2$ be a graph of order $n  \ (n\geq 5)$ obtained from a cycle of order $4$ and four paths
$P_1,P_2,P_3,P_4$ of length $a,b,c,d \ (0\leq a,b,c,d\leq n-4, \ a+b+c+d=n-4)$ respectively by identifying
each vertex of this cycle with an endvertex of one of the four
paths.

\begin{figure}[!hbpt]
\begin{center}
\includegraphics[scale=0.75]{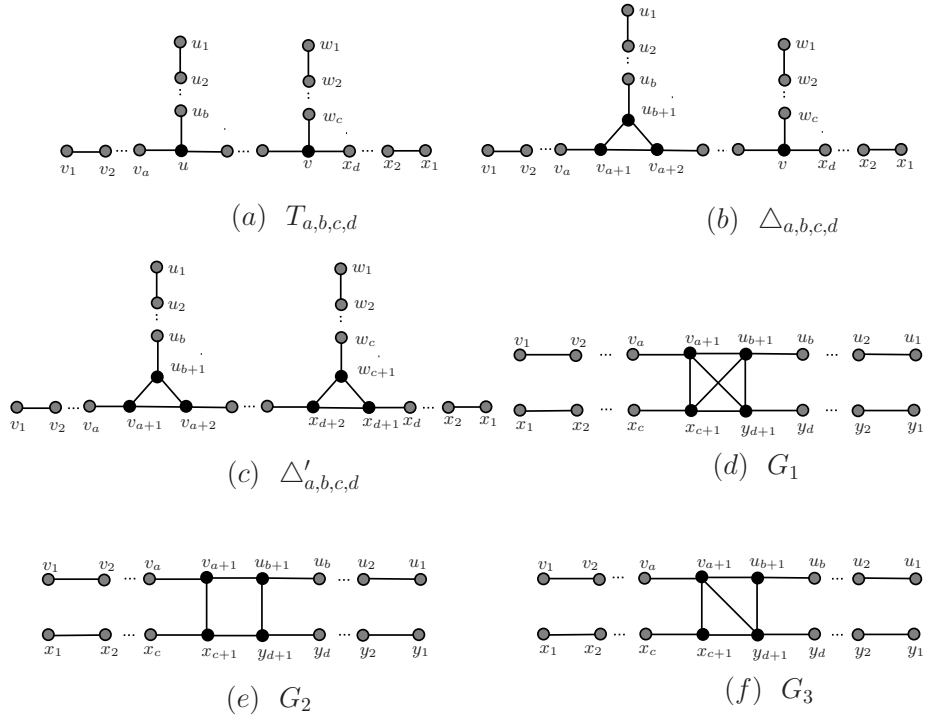}\\
\caption{Graphs for Theorem \ref{th4}.}
\end{center}\label{fig2}
\end{figure}

\item Let $G_3$ be a graph of order $n  \ (n\geq 5)$ obtained from $K_4^-$ and four paths
$P_1,P_2,P_3,P_4$ of length $a,b,c,d \ (0\leq a,b,c,d\leq n-4, \ a+b+c+d=n-4)$ respectively by identifying
each vertex of $K_4^-$ with an endvertex of one of the four paths,
where $K_4^-$ denotes the graph obtained from a clique of order $4$
by deleting one edge.
\end{itemize}

\begin{thm}\label{th4}
Let $G$ be a connected graph of order $n \ (n\geq 5)$. Then $sdiam_4(G)=n-1$ if
and only if $G=T_{a,b,c,d}$ or $G=\triangle_{a,b,c,d}$ or
$G=\triangle^\prime_{a,b,c,d}$ or $G=G_1$ or $G=G_2$ or $G=G_3$.
\end{thm}

\section{Proofs of Theorem $2$ and $3$}

In this section, we characterize graphs with $sdiam_4(G)=3,4$ and give the proofs of Theorems $2$ and $3$.
\begin{lem}\label{lem1}
Let $G$ be a connected graph of order $n$, and let $k$ be an integer with $3\leq k\leq n-1$. Then $sdiam_k(G)=n-1$ if and only if the number of non-cut vertices in
$G$ is at most $k$.
\end{lem}
\begin{pf}
Let $r$ be the number of non-cut vertices in $G$.
Suppose $sdiam_k(G)=n-1$. We claim that $r\leq k$. Assume, to the
contrary, that $r\geq k+1$. For any $S\subseteq V(G)$ with $|S|=k$,
there exists a non-cut vertex in $G$, say $u$, such that $u\in
V(G)\setminus S$. Then $G\setminus u$ is connected, and hence
$G\setminus u$ contains a spanning tree of size $n-2$. From the
arbitrariness of $S$, we have $sdiam_{k}(G)\leq d_{T}(S)\leq n-2$, a
contradiction. So $r\leq k$, as desired.

Conversely, we suppose $r\leq k$. Let $v_1,v_2,\cdots,v_r$ be all
the non-cut vertices in $G$. Then the remaining vertices are all cut
vertices of $G$. Choose $v_{i_1},v_{i_2},\cdots,v_{i_{k-r}}\in V(G)\setminus
\{v_1,v_2,\cdots,v_r\}$. Set
$S=\{v_1,v_2,\cdots,v_r,v_{i_1},v_{i_2},\cdots,v_{i_{k-r}}\}$. Note that each
vertex in $V(G)\setminus S$ is a cut vertex of $G$. Therefore, any
$S$-Steiner tree $T$ occupies all the vertices of $G$, and hence $sdiam_k(G)\geq d_G(S)\geq n-1$. From Theorem \ref{th1}, we
have $sdiam_k(G)=n-1$, as desired.\qed
\end{pf}

\vskip 0.3cm

The following corollary is immediate from the above lemma.
\begin{cor}\label{cor1}
Let $G$ be a connected graph of order $n$, and let $k$ be an integer with $3\leq k\leq n-2$. Then $sdiam_k(G)\leq n-2$ if
and only if the number of non-cut vertices in
$G$ is at least $k+1$.
\end{cor}

Mao \cite{Mao} obtained the following result, which will be used
later.

\begin{lem}{\upshape \cite{Mao}}\label{lem2}
Let $n,k$ be two integers with $2\leq k\leq n$, and let $G$ be a
connected graph of order $n$. If $sdiam_k(G)=k-1$, then $0\leq
\Delta(\overline{G})\leq k-2$, namely, $n-k+1\leq \delta(G)\leq
n-1$.
\end{lem}

\noindent \textbf{Proof of Theorem \ref{th2}:} If $n=4$, then $sdiam_4(G)=3$. So we assume that $n\geq 5$.
Suppose
$sdiam_4(G)=3$. For Lemma \ref{lem2}, if $sdiam_4(G)=3$, then
$n-3\leq\delta(G)\leq n-1$. We claim that $C_4$ is not a subgraph of $\overline{G}$. Assume, to the contrary, that
$C_4$ is a subgraph of $\overline{G}$. Choose $S=V(C_4)$.
Since $G[S]$ is not connected, it follows that any $S$-Steiner tree
must contain one vertex in $V(G)\setminus S$, and hence
$sdiam_4(G)\geq d_G(S)\geq 4$, a contradiction. So $C_4$ is not a subgraph of $\overline{G}$.

Conversely, we suppose that $n-3\leq\delta(G)\leq n-1$ and
$C_4$ is not a subgraph of $\overline{G}$. Since
$n-3\leq\delta(G)\leq n-1$, it follows that $G$ is a graph obtained
from the complete graph of order $n$ by deleting some independent
paths and cycles. For any $S\subseteq V(G)$, since
$C_4$ is not a subgraph of $\overline{G}$, it follows
that $\overline{G}[S]=4K_1$ or $\overline{G}[S]=K_2\cup 2K_1$ or
$\overline{G}[S]=2K_2$ or $\overline{G}[S]=P_4$ or
$\overline{G}[S]=K_3\cup K_1$ or $\overline{G}[S]=P_3\cup K_1$. Then
$G[S]=K_4$ or $G[S]=K_4\setminus e$ or $G[S]=C_4$ or $G[S]=P_4$ or
$G[S]=K_{1,3}$ or $G[S]=K_{1,3}^{+}$, where $K_{1,3}^{+}$ is the
graph obtained from a star $K_{1,3}$ by adding an edge. Since $G[S]$ is a connected graph,
it follows that $d_G(S)\leq 3$. From the arbitrariness of $S$, we
have $sdiam_4(G)\leq 3$ and hence $sdiam_4(G)=3$ by Theorem
\ref{th1}. The proof is complete.\qed \\

\noindent \textbf{Proof of Theorem \ref{th3}:} Suppose that $G$ is a
graph with $sdiam_4(G)=4$. From Theorem \ref{th2}, we have $\delta(G)=n-3$ and $C_4$ is a subgraph of $\overline{G}$, or $\delta(G)\leq n-4$.
For the former, we have $\delta(G)=n-3$ and $C_4$ is a subgraph of $\overline{G}$, as desired. Suppose $\delta(G)\leq n-4$. It suffices to prove that each $H_i\ (1\leq i\leq 4)$ is not a spanning subgraph of $\overline{G}$, and we have the following claims.

\textbf{Claim 1.} $H_1$ is not a spanning subgraph of $\overline{G}$.

\noindent\textbf{Proof of Claim 1.} Assume, to the contrary, that
$H_1$ is a spanning subgraph of $\overline{G}$. Choose
$S=\{u_1,u_2,u_3,u_4\}\subseteq V(H_1)=V(\overline{G})$. Then the subgraph in $\overline{G}$ induced
by the vertices in $S$ is a complete graph of order $4$, and hence $G[S]=4K_1$ is not connected. Therefore, any $S$-Steiner tree $T$
must occupy a vertex in $V(G)\setminus S$, say $x$. Because
$H_1$ is a spanning subgraph of $\overline{G}$, we have $xu_1\notin
E(G)$ or $xu_2\notin E(G)$ or $xu_3\notin E(G)$ or $xu_4\notin
E(G)$. Thus, the $S$-Steiner tree $T$ must occupy another vertex in
$V(G)\setminus S$, and hence the tree $T$ must occupy at
least two vertices in $V(G)\setminus S$. Then $d_G(S)\geq 5$,
and hence $sdiam_4(G)\geq 5$, a contradiction. So $H_1$ is not a spanning subgraph of $\overline{G}$, as desired.\qed

\textbf{Claim 2.} $H_2$ is not a spanning subgraph of $\overline{G}$.

\noindent\textbf{Proof of Claim 2.} Assume, to the contrary, that
$H_2$ is a spanning subgraph of $\overline{G}$. Choose
$S=\{u_1,u_2,u_3,u_4\}\subseteq V(H_2)=V(\overline{G})$. Since $G[S]$ is not connected, it follows
that any $S$-Steiner tree $T$ must occupy a vertex in $V(G)\setminus
S$, say $x$. From the structure of $H_2$, since $H_2$ is a spanning subgraph of $\overline{G}$, we have $xu_1,xu_4\in
E(\overline{G})$ or $xu_2\in E(\overline{G})$ or $xu_3\in
E(\overline{G})$. If $xu_1,xu_4\in E(\overline{G})$, then there are
at most three edges in $\{xu_2,xu_3,u_1u_4\}$ belonging to $G[S\cup
\{x\}]$. In order to connect to $u_1$ or $u_4$, the $S$-Steiner tree
$T$ uses at least two vertex of $V(G)\setminus S$. If $xu_2\in
E(\overline{G})$, then there are at most three edges in
$\{xu_1,xu_3,xu_4,u_1u_4\}$ belonging to $G[S\cup \{x\}]$. In order
to connect to $u_2$, the $S$-Steiner tree $T$ must use at least
two vertex of $V(G)\setminus S$. The same is true for $xu_3\in
E(\overline{G})$. Therefore, $e(T)\geq 5$ and $d_G(S)\geq 5$, which
results in $sdiam_4(G)\geq 5$, a contradiction. So $H_2$ is not a spanning subgraph of $\overline{G}$.\qed

\textbf{Claim 3.} $H_3$ is not a spanning subgraph of $\overline{G}$.

\noindent\textbf{Proof of Claim 3.} Assume, to the contrary, that
$H_3$ is a spanning subgraph of $\overline{G}$. Choose
$S=\{u_1,u_2,u_3,u_4\}\subseteq V(H_3)=V(\overline{G})$. Since $G[S]=2K_2$ is not connected, it
follows that any $S$-Steiner tree $T$ must occupy a vertex in
$V(G)\setminus S$, say $x$. From the structure of $H_3$, since
$H_3$ is a spanning subgraph of $\overline{G}$, we have $xu_1,xu_2\in
E(\overline{G})$ or $xu_3,xu_4\in E(\overline{G})$. If $xu_1,xu_2\in
E(\overline{G})$, then there are at most four edges in
$\{xu_3,xu_4,u_1u_2,u_3u_4\}$ belonging to $G[S\cup \{x\}]$. In
order to connect to $u_1$ or $u_2$, the $S$-Steiner tree $T$ uses at
least two vertex of $V(G)\setminus S$. If $xu_3,xu_4\in
E(\overline{G})$, then there are at most four edges in
$\{xu_1,xu_2,u_1u_2,u_3u_4\}$ belonging to $G[S\cup \{x\}]$. In
order to connect to $u_3$ or $u_4$, the $S$-Steiner tree $T$ must use at
least two vertex of $V(G)\setminus S$. Therefore, $e(T)\geq 5$ and
$d_G(S)\geq 5$, which results in $sdiam_4(G)\geq 5$, a
contradiction. So $H_3$ is not a spanning subgraph of $\overline{G}$.\qed

\textbf{Claim 4.} $H_4$ is not a spanning subgraph of $\overline{G}$.

\noindent\textbf{Proof of Claim 4.} Assume, to the contrary, that
$H_4$ is a spanning subgraph of $\overline{G}$. Choose
$S=\{u_1,u_2,u_3,u_4\}\subseteq V(H_4)=V(\overline{G})$. Since $G[S]=K_3\cup K_1$ is not connected,
it follows that any $S$-Steiner tree $T$ must occupy a vertex in
$V(G)\setminus S$, say $x$. From the structure of $H_4$, since
$H_4$ is a spanning subgraph of $\overline{G}$, we have $xu_3\in
E(\overline{G})$ or $xu_1,xu_2,xu_4\in E(\overline{G})$. If $xu_3\in
E(\overline{G})$, then there are at most six edges in
$\{xu_1,xu_2,xu_4,u_1u_2,u_1u_4,u_2u_4\}$ belonging to $G[S\cup
\{x\}]$. In order to connect to $u_3$, the $S$-Steiner tree $T$ uses
at least two vertex of $V(G)\setminus S$. If $xu_1,xu_2,xu_4\in
E(\overline{G})$, then there are at most four edges in
$\{xu_3,u_1u_2,u_1u_4,u_2u_4\}$ belonging to $G[S\cup \{x\}]$. In
order to connect to $u_1$ or $u_2$ or $u_4$, the $S$-Steiner tree
$T$ uses at least two vertex of $V(G)\setminus S$. Therefore,
$e(T)\geq 5$ and $d_G(S)\geq 5$, which results in $sdiam_4(G)\geq
5$, a contradiction. So $H_4$ is not a spanning subgraph of $\overline{G}$.\qed

From the above argument, we know that the result holds.

Conversely, suppose that $G$ is a connected graph satisfying one of
the following conditions.

$\bullet$ $\delta(G)=n-3$ and $C_4$ is a subgraph of $\overline{G}$;

$\bullet$ $\delta(G)\leq n-4$ and $H_i \ (1\leq
i\leq 4)$ is not a spanning subgraph of $\overline{G}$.

Suppose that $\delta(G)=n-3$ and $C_4$ is a subgraph of $\overline{G}$. Since $\delta(G)=n-3$, it follows that $G$ is a graph
obtained from the complete graph of order $n$ by deleting some
pairwise independent paths and cycles. Then $\overline{G}$ is a union of pairwise independent paths, cycles, and isolated vertices. For any $S=\{u,v,w,z\}\subseteq V(G)$,
since $\overline{G}$ contains $C_4$ as its subgraph, it follows that
$\overline{G}[S]=C_4$ or $\overline{G}[S]=4K_1$ or
$\overline{G}[S]=K_2\cup 2K_1$ or $\overline{G}[S]=2K_2$ or
$\overline{G}[S]=P_4$ or $\overline{G}[S]=K_3\cup K_1$ or
$\overline{G}[S]=P_3\cup K_1$. Then $G[S]=2K_2$ or $G[S]=K_4$ or
$G[S]=K_4\setminus e$ or $G[S]=C_4$ or $G[S]=P_4$ or $G[S]=K_{1,3}$ or
$G[S]=K_{1,3}^{+}$, where $K_{1,3}^{+}$ is the graph obtained from $K_{1,3}$ by
adding an edge. If $G[S]=2K_2$, then $|E_G[x,S]|=4$ for any $x$ in
$V(G)\setminus S$, since $\delta(G)=n-3$. Thus, we have $d_G(S)=4$. For the other cases, $G[S]$
is connected, and so $d_G(S)=3$. From the arbitrariness of $S$, we
have $sdiam_4(G)=4$, as desired.

Suppose $\delta(G)\leq n-4$ and each $H_i \ (1\leq
i\leq 4)$ is not a spanning subgraph of $\overline{G}$.
For any $S\subseteq V(G)$ and $|S|=4$, if there exists a vertex $x\in V(G)\setminus S$ such that $|E(\overline{G}[x,S])|=0$,
then $|E(G[x,S])|=4$, and hence
the tree $T$ induced by the four edges in $E(G[x,S])$ is an $S$-Steiner
tree in $G$, and hence $d_G(S)\leq 4$, as desired. From now on, we assume for any $S\subseteq V(G)$ and $|S|=4$,
and any $x\in V(G)\setminus S$, $|E(\overline{G}[x,S])|\geq
1$.

From the definition of $sdiam_4(G)$ and
Theorem \ref{th2}, it suffices to show that $d_G(S)\leq 4$ for any
set $S\subseteq V(G)$ and $|S|=4$. It is clear that $0\leq
|E(G[S])|\leq 6$. If $4\leq |E(G[S])|\leq 6$, then $G[S]$ is
connected, and hence $G[S]$ contains a spanning tree, which is an
$S$-Steiner tree in $G$. So $d_G(S)=3<4$, as desired. From now on,
we assume $0\leq |E(G[S])|\leq 3$.

If $|E(G[S])|=0$, then $\overline{G}[S]=K_4$. Since
$|E(\overline{G}[x,S])|\geq 1$ for any $x\in V(G)\setminus S$, it
follows that $H_1$ is a spanning subgraph of $\overline{G}$, a contradiction.

Suppose $|E(G[S])|=1$. Set $S=\{u_1,u_2,u_3,u_4\}$. Without loss of
generality, let $u_1u_4\in E(G)$ and
$u_1u_2,u_1u_3,u_2u_3,u_2u_4,u_3u_4\notin E(G)$. Then
$u_1u_2,u_1u_3,u_2u_3,u_2u_4,u_3u_4\in E(\overline{G})$, and hence $\overline{G}[S]$ is a graph
obtained from $K_4$ by deleting one edge. Since $H_2$ is not a spanning subgraph of $\overline{G}$,
it follows that there exists a vertex $x\in
V(G)-S$ such that $xu_1\in E(\overline{G})$ but
$xu_2,xu_3,xu_4\notin E(\overline{G})$ or $xu_4\in E(\overline{G})$
but $xu_2,xu_3,xu_1\notin E(\overline{G})$. By symmetry, we only to
consider the former case. Clearly, $xu_2,xu_3,xu_4\in E(G)$. Combining
this with $u_1u_4\in E(G)$, the tree $T$ induced by the edges in
$\{u_1u_4,xu_2,xu_3,xu_4\}$ is an $S$-Steiner tree in $G$ and hence
$d_G(S)\leq 4$, as desired.

Suppose $|E(G[S])|=2$. Without loss of generality, we can assume that
$u_1u_2,u_3u_4\in E(G)$ or $u_1u_2,u_1u_4\in E(G)$. First, we
consider the case $u_1u_2,u_3u_4\in E(G)$. Clearly,
$u_1u_3,u_1u_4,u_2u_3,u_2u_4\notin E(G)$, and hence
$u_1u_3,u_1u_4,u_2u_3,u_2u_4\in E(\overline{G})$. Note that for any $S\subseteq V(G)$ and $|S|=4$,
and any $x\in V(G)\setminus S$, $|E(\overline{G}[x,S])|\geq
1$. Since $H_3$ is not a spanning subgraph of $\overline{G}$, it follows that
there exists a vertex $x\in V(G)\setminus S$ satisfying one of the
following.
\begin{itemize}
\item[] (1) $xu_1\in E(\overline{G})$ but $xu_2,xu_3,xu_4\notin
E(\overline{G})$;

\item[] (2) $xu_1,xu_4\in E(\overline{G})$ but
$xu_2,xu_3\notin E(\overline{G})$;

\item[] (3) $xu_1,xu_3\in E(\overline{G})$ but $xu_2,xu_4\notin
E(\overline{G})$;

\item[] (4) $xu_2\in E(\overline{G})$ but
$xu_1,xu_3,xu_4\notin E(\overline{G})$;

\item[] (5) $xu_2,xu_4\in E(\overline{G})$ but
$xu_1,xu_3\notin E(\overline{G})$;

\item[] (6) $xu_2,xu_3\in
E(\overline{G})$ but $xu_1,xu_4\notin E(\overline{G})$;

\item[] (7) $xu_3\in E(\overline{G})$
but $xu_1,xu_2,xu_4\notin E(\overline{G})$;

\item[] (8) $xu_3,xu_1\in
E(\overline{G})$ but $xu_2,xu_4\notin E(\overline{G})$;

\item[] (9) $xu_2,xu_3\in E(\overline{G})$ but $xu_1,xu_4\notin
E(\overline{G})$;

\item[] (10) $xu_4\in E(\overline{G})$ but
$xu_1,xu_2,xu_3\notin E(\overline{G})$;

\item[] (11) $xu_1,xu_4\in E(\overline{G})$ but $xu_2,xu_3\notin
E(\overline{G})$;

\item[] (12) $xu_2,xu_4\in E(\overline{G})$ but
$xu_1,xu_3\notin E(\overline{G})$.
\end{itemize}

By symmetry, we only consider the first three cases, and other cases can be similarly proved. If $xu_1\notin E(G)$ but $xu_2,xu_3,xu_4\in E(G)$, then the tree $T$ induced by the edges in $\{xu_2,xu_4,u_1u_2,u_3u_4\}$ is
an $S$-Steiner tree in $G$ and hence $d_G(S)\leq 4$. If $xu_1,xu_4\notin
E(G)$ but $xu_2,xu_3\in E(G)$, then the tree $T$ induced by the edges in $\{u_1u_2,u_3u_4,xu_2,xu_3\}$ is
an $S$-Steiner tree in $G$ and hence $d_G(S)\leq 4$. If $xu_1,xu_3\notin E(G)$ but
$xu_2,xu_4\in E(G)$, then the tree $T$ induced by the edges in $\{u_1u_2,u_3u_4,xu_2,xu_4\}$ is
an $S$-Steiner tree in $G$ and hence $d_G(S)\leq 4$, as desired.

Next, we consider the case $u_1u_2,u_1u_4\in E(G)$. Clearly,
$u_1u_3,u_2u_3,u_2u_4,u_3u_4\notin E(G)$, and hence
$u_1u_3,u_2u_3,u_2u_4,u_3u_4\in E(\overline{G})$. Note that for any $S\subseteq V(G)$ and $|S|=4$,
and any $x\in V(G)\setminus S$, $|E(\overline{G}[x,S])|\geq
1$. Since $H_4$ is not a spanning subgraph of $\overline{G}$, it follows that there exists a vertex $x\in
V(G)\setminus S$ satisfying one of the following.
\begin{itemize}
\item[] (1) $xu_1\in E(\overline{G})$ but $xu_2,xu_3,xu_4\notin
E(\overline{G})$;

\item[] (2) $xu_1,xu_4\in E(\overline{G})$ but
$xu_2,xu_3\notin E(\overline{G})$;

\item[] (3) $xu_1,xu_2\in E(\overline{G})$ but $xu_4,xu_3\notin
E(\overline{G})$;

\item[] (4) $xu_2\in E(\overline{G})$ but
$xu_1,xu_3,xu_4\notin E(\overline{G})$;

\item[] (5) $xu_2,xu_4\in E(\overline{G})$ but $xu_1,xu_3\notin
E(\overline{G})$;

\item[] (6) $xu_1,xu_2\in E(\overline{G})$ but
$xu_4,xu_3\notin E(\overline{G})$;

\item[] (7) $xu_4\in E(\overline{G})$ but
$xu_1,xu_2,xu_3\notin E(\overline{G})$;

\item[] (8) $xu_1,xu_4\in
E(\overline{G})$ but $xu_2,xu_3\notin E(\overline{G})$;

\item[] (9) $xu_2,xu_4\in E(\overline{G})$ but
$xu_1,xu_3\notin E(\overline{G})$.
\end{itemize}

By symmetry, we only consider the first three cases. Note that $u_1u_2,u_1u_4\in E(G)$.
If $xu_1\notin E(G)$ but $xu_2,xu_3,xu_4\in E(G)$, then the tree $T$ induced by the edges in $\{u_1u_2,u_1u_4,xu_2,xu_3\}$ is
an $S$-Steiner tree in $G$ and hence $d_G(S)\leq 4$.
If $xu_1,xu_4\notin
E(G)$ but $xu_2,xu_3\in E(G)$, then the tree $T$ induced by the edges in $\{u_1u_2,u_1u_4,xu_2,xu_3\}$ is
an $S$-Steiner tree in $G$ and hence $d_G(S)\leq 4$.
If $xu_1,xu_2\notin E(G)$ but
$xu_4,xu_3\in E(G)$, then the tree $T$ induced by the edges in $\{u_1u_2,u_1u_4,xu_4,xu_3\}$ is
an $S$-Steiner tree in $G$ and hence $d_G(S)\leq 4$, as desired.

Suppose $|E(G[S])|=3$. Without loss of generality, let
$u_1u_2,u_1u_4,u_2u_3\in E(G)$ or $u_1u_2,u_1u_4,u_2u_4\in E(G)$. If
$u_1u_2,u_1u_4,u_2u_3\in E(G)$, then the tree $T$ induced by the
edges in $\{u_1u_2,u_1u_4,u_2u_3\}$ is an $S$-Steiner tree in $G$ and
hence $d_G(S)=3\leq 4$, as desired. If $u_1u_2,u_1u_4,u_2u_4\in
E(G)$, then $u_1u_3,u_2u_4,u_2u_3\notin E(G)$, and hence
$u_1u_3,u_2u_4,u_2u_3\in E(\overline{G})$. Since $H_4$ is not a spanning subgraph of $\overline{G}$, it follows that there exists a vertex $x\in
V(G)\setminus S$ satisfying one of the following.
\begin{itemize}
\item[] (1) $xu_1\in E(\overline{G})$ but $xu_2,xu_3,xu_4\notin
E(\overline{G})$;

\item[] (2) $xu_1,xu_4\in E(\overline{G})$ but
$xu_2,xu_3\notin E(\overline{G})$;

\item[] (3) $xu_1,xu_2\in E(\overline{G})$ but $xu_3,xu_4\notin
E(\overline{G})$;

\item[] (4) $xu_2\in E(\overline{G})$ but
$xu_1,xu_3,xu_4\notin E(\overline{G})$;

\item[] (5) $xu_2,xu_4\in E(\overline{G})$ but $xu_1,xu_3\notin
E(\overline{G})$;

\item[] (6) $xu_1,xu_2\in E(\overline{G})$ but
$xu_3,xu_4\notin E(\overline{G})$;

\item[] (7) $xu_4\in E(\overline{G})$ but
$xu_2,xu_3,xu_1\notin E(\overline{G})$;

\item[] (8) $xu_1,xu_4\in
E(\overline{G})$ but $xu_2,xu_3\notin E(\overline{G})$;

\item[] (9) $xu_4,xu_2\in E(\overline{G})$ but
$xu_3,xu_1\notin E(\overline{G})$.
\end{itemize}

By symmetry, we only consider the first three cases. Recall that $u_1u_2,u_1u_4,u_2u_4\in
E(G)$. If $xu_1\notin E(G)$ but $xu_2,xu_3,xu_4\in E(G)$, then the tree $T$ induced by the edges in $\{u_1u_2,u_1u_4,xu_2,xu_3\}$ is an $S$-Steiner tree in $G$ and hence $d_G(S)\leq 4$. If $xu_1,xu_4\notin
E(G)$ but $xu_2,xu_3\in E(G)$, then the tree $T$ induced by the edges in $\{u_1u_2,u_1u_4,xu_2,xu_3\}$ is
an $S$-Steiner tree in $G$ and hence $d_G(S)\leq 4$. If $xu_1,xu_2\notin E(G)$ but
 $xu_3,xu_4\in E(G)$, then the tree $T$ induced by the edges in $\{u_1u_2,u_1u_4,xu_3,xu_4\}$ is
an $S$-Steiner tree in $G$ and hence $d_G(S)\leq 4$.

From the arbitrariness of $S$, we have $sdiam_4(G)\leq 4$. Since
$\delta(G)=n-3$ and $C_4\in\overline{G}$, or $\delta(G)\leq n-4$, it follows from
Theorem \ref{th2} that $sdiam_4(G)=4$. The proof is now
complete.\qed

\section{Proof of Theorem $4$}

The following lemma is a preparation of our main result.
\begin{pro}\label{pro1}
Let $G$ be a connected graph, and $H$ be a connected subgraph
of $G$. Then the number of non-cut vertices of $G$ is not less than the
number of non-cut vertices of $H$.
\end{pro}
\begin{pf}
It suffices to show that there exists an injective mapping $f$ from
the set of non-cut vertices of $H$ to the set of non-cut vertices of $G$. We
define such a mapping $f$ as follows. Let $v$ be a non-cut vertex of $H$.
If $v$ is a non-cut vertex of $G$, then let $f(v) = v$. If $v$ is a cut-vertex
of $G$, then let $G_1$ be a component of $G\setminus v$ not containing any vertex
of $H$. Let $T_1$ be a spanning tree of $G_1$, and let $w$ be an end-vertex of
$T_1$ distinct from $v$. Then $w$ is a non-cut vertex of $G$, and we define
$f(v)=w$. Now $f$ maps non-cut vertices of $H$ to non-cut vertices of $G$,
and $f$ is injective since either $f(v)=v$ or $f(v$) is in a component of
$G\setminus V(H)$ which is (in $G$) attached only to $v$, and to no other vertex
in $V(H)$.\qed
\end{pf}

\vskip 0.3cm

From Proposition \ref{pro1}, the following corollaries are immediate.
\begin{cor}\label{cor2}
Let $G$ be a connected graph of order $n \ (n\geq 3)$, and let $c(G)$ be the
circumference of the graph $G$. If $c(G)\leq n$, then there are at
least $c(G)$ non-cut vertices in $G$.
\end{cor}

\begin{cor}\label{cor3}
Let $G$ be a connected graph of order $n  \ (n\geq 3)$. Let $C_1,C_2,\cdots,C_r \ (r\geq 2)$
are cycles of the graph $G$ with $|V(C_i)|=n_i \ (1\leq i\leq r)$.
If $|V(C_i)\cap V(C_j)|\leq 1$ for any $i,j \ (1\leq i,j\leq r, \
i\neq j)$, then the graph $G$ has at least
$n_1+n_2+\cdots+n_r-2(r-1)$ are non-cut vertices in $G$.
\end{cor}

We are now in a position to give the proof of Theorem \ref{th4}.

\noindent \textbf{Proof of Theorem \ref{th4}:} Suppose
$G=T_{a,b,c,d}$ or $G=\triangle_{a,b,c,d}$ or
$G=\triangle^\prime_{a,b,c,d}$ or $G=G_1$ or $G=G_2$ or $G=G_3$.
Since there are at most four non-cut vertices in $G$, it follows
from Lemma \ref{lem1} that $sdiam_4(G)=n-1$.

Conversely, suppose $sdiam_4(G)=n-1$. If $G$ is a tree, then it
follows from Lemma \ref{lem1} that $G$ contains at most non-cut four
vertices, and hence $G=T_{a,b,c,d}$. Now, we assume that $G$ contains
cycles. Recall that $c(G)$ is the circumference of the graph $G$.
Obviously, $3\leq c(G)\leq n$. If $5\leq c(G)\leq n$, then it follow
from Corollaries \ref{cor1} and \ref{cor2} that $sdiam_4(G)\leq
n-2$, a contradiction. Therefore, $c(G)=3$ or $c(G)=4$. If $c(G)=4$,
then it follows from Lemma \ref{lem1} and Corollaries \ref{cor1} and \ref{cor3} that $G$ contains four non-cut vertices, and from Corollaries \ref{cor3} that $G$
contains no two cycles $C_1,C_2$ with $|V(C_1)|=4$ or $|V(C_2)|=4$ such that $|V(C_1)\cap V(C_2)|\leq
1$.
From Proposition \ref{pro1}, we have the following facts.

$\bullet$ $G\setminus V(C_i) \ (i=1,2)$ is a union of pairwise independent paths;

$\bullet$ The number of these paths are at most four;

$\bullet$ The endvertices of each pair of these paths share the different neighbors in $C_i$.

From these facts, we have $G=G_1$ or $G=G_2$ or $G=G_3$. If $c(G)=3$, then it
follows from Lemma \ref{lem1} and Corollaries
\ref{cor1} and \ref{cor3} that $G$ contains at most four non-cut vertices, and $G$ contains exactly one triangle or at most two cycles $C_1$ and $C_2$ with $|V(C_1)|=3$ and $|V(C_2)|=3$ such that
$|V(C_1)\cap V(C_2)|\leq 1$. Suppose $G$ contains only one triangle. Let $K_{1,3}^*$ be the subdivision of star $K_{1,3}$ of order $t$. Then we have the following facts.

$\bullet$ The graph obtained from $G$ by deleting this triangle is $P_r\cup P_s\cup K_{1,3}^* \ (r+s+t=n-3)$ or $P_r\cup K_{1,3}^* \ (r+t=n-3)$ or $P_r\cup P_s \ (r+s=n-3)$ or $P_{n-3}$ or $K_{1,3}^* \ (t=n-3)$ or $P_r\cup P_s\cup P_p\ (r+s+p+q=n-3)$ or $P_r\cup P_s\cup P_p\cup P_q\ (r+s+p+q=n-3)$;

$\bullet$ If the graph obtained from $G$ by deleting this triangle is not $P_r\cup P_s\cup P_p\cup P_q\ (r+s+p+q=n-3)$, then the endvertices of each pair of these paths share the different neighbors in the triangle.

$\bullet$ If the graph obtained from $G$ by deleting this triangle is $P_r\cup P_s\cup P_p\ (r+s+p+q=n-3)$ or $P_r\cup P_s\cup P_p\cup P_q\ (r+s+p+q=n-3)$, then each vertex of this triangle share at least one common neighbor of each path.

Clearly, we have $G=\triangle_{a,b,c,d}$. If $G$ contains at most two cycles $C_1$ and $C_2$ with $|V(C_1)|=3$ and $|V(C_2)|=3$ such that $|V(C_1)\cap V(C_2)|\leq 1$, then $G\setminus (V(C_1)\cup V(C_2)$ is a union of pairwise independent paths, and the number of these path is at most five. So $G=\triangle^\prime_{a,b,c,d}$. The proof is complete.\qed

\end{document}